\documentclass[12pt]{article}
 \usepackage{amssymb}
 \usepackage{amsmath}
 \usepackage[latin1]{inputenc}
 \usepackage[T1]{fontenc}
\usepackage{times}
 \usepackage{dsfont}
 \usepackage{amsfonts}
 \usepackage{natbib}
 \usepackage{mathrsfs}
\date{}

 \usepackage{graphicx}

\newtheorem{theorem}{Theorem}[section]

\newtheorem{e-proposition}[theorem]{Proposition}

\newtheorem{e-definition}[theorem]{Definition\rm}
\newtheorem{remark}{\it Remark\/}

\begin{document}
 \begin{center}
 {\LARGE {Large Deviations Theorems in Nonparametric Regression on Functional Data}}
 \end{center}
 \begin{center}
 \large{Mohamed Cherfi}\vspace{3mm}

 {\small \it Laboratoire de Statistique Th\'eorique et Appliqu\'ee (LSTA)\\ Equipe d'Accueil 3124\\ Universit\'e Pierre et Marie Curie -- Paris 6\\ Tour 15-25, 2\`eme \'etage\\ 4 place Jussieu\\  75252 Paris cedex 05 }
 \end{center}
\renewcommand{\thefootnote}{}
\footnote{ \vspace*{-4mm} \noindent{{\it E-mail address:}
mohamed.cherfi@gmail.com} }

 \vspace{3mm} \hrule \vspace{3mm} {\small \noindent{\bf Abstract.}
\\
In this Note we prove large deviations principles for the
Nadaraya-Watson estimator of the regression of a real-valued variable with a functional covariate. Under suitable conditions, we show pointwise and uniform large deviations theorems with  good rate functions.\\

\vspace{0.2cm}
 
 \vspace{4mm}\hrule
\section{Introduction}
Let $\{(Y_i,\,X_i),~i\geq 1\}$ be a sequence of independent and identically distributed
random vectors. The random variables $Y_i$ are real, with $\mathds{E}|Y|<\infty$, and the $X_i$ are random vectors with values in a semi-metric space $(\mathcal{X},d(\cdot,\cdot))$.

Consider now the functional regression model,
\begin{equation}\label{eqReg}
   Y_i:=\mathds{E}(Y|X_i)+\varepsilon_i=r(X_i)+\varepsilon_i \quad i=1,\ldots,n,
\end{equation}
where $r$ is the regression operator mapping $\mathcal{X}$ onto $\mathbb{R}$, and the $\varepsilon_i$ are real variables such that, for all
$i$, $\mathds{E}(\varepsilon_i|X_i ) = 0$ and $\mathds{E}(\varepsilon_i^2|X_i ) =\sigma^2_{\varepsilon}(X_i ) < \infty$. Note that in practice $\mathcal{X}$ is a normed
space which can be of infinite dimension (e.g., Hilbert or Banach space) with norm $\|\cdot\|$ so that
$d(x, x') = \|x - x'\|$, which is the case in this paper.

\noindent \cite{FerratyVieu2004} provided a consistent estimate for the nonlinear regression operator $r$, based on the usual finite-dimensional smoothing ideas,  that is
\begin{equation}\label{eqNW}
   \displaystyle{\widehat{r}_n(x):=\frac{\sum_{i=1}^n Y_iK\bigg(\frac{\|x-X_i\|}{h_n}\bigg)}{\sum_{i=1}^nK\bigg(\frac{\|X_i-x\|}{h_n}\bigg)}},
\end{equation}
where $K(\cdot)$ is a real-valued kernel and $h_n$ the bandwidth, is a sequence of positive real numbers converging to $0$ as $n\longrightarrow\infty$. Note that the bandwidth $h_n$ depends on $n$, but we drop this index for notational convenience. In what follows $K_{h}(u)$ stands for $\displaystyle{K\bigg(\frac{u}{h}\bigg)}$. The estimator defined in (\ref{eqNW}) is a generalization to the functional framework of the classical Nadaraya-Watson regression estimator. The asymptotic properties of this estimate have been studied extensively by several authors, we cite among others \cite{FerratyMasVieu2007}, for a complete survey see the monograph by \cite{FerratyVieu2006}.

\noindent There exists a wide large deviations literature, for an extensive overview on the area
we refer to the monograph of \cite{DemboZeitouni1998}. The large deviations behavior of the Nadaraya-Watson estimate of the regression function, have been studied at first by  \cite{Louani1999}, sharp results have been obtained by \cite{Joutard2006} in the univariate framework. In the multidimensional case \cite{MokkademPelletierThiam2008} obtained pointwise large and moderate deviations results for the Nadaraya-Watson and recursive kernel estimators of  the regression. 

\noindent Large deviations results have many implications in statistical analysis, we refer  to \cite{Bahadur1971} for a good survey on the subject. The main concern in the nonparametric setup is with asymptotic efficiency and the inaccuracy rate. These questions have been investigated in \cite{Louani1998} for the density case. As an application  of large deviations results, \cite{Louani1999} gives the inaccuracy rate in conditional distribution function estimation. Another application can be found in \cite{BerrahouLouani2006}, where the authors provide Bahadur efficiency of symmetry tests. We refer to the monograph of \cite{Nikitin1995} for an account of results on the asymptotic efficiency.

\noindent All the results of large deviations in nonparametric regression only dealt with the case of scalar covariates. In this Note, our aim is to obtain such results in the functional case. Inaccuracy rate in nonparametric regression estimation on functional data will be considered in a futur work. We are also interested to obtain a minimax Bahadur-type risk of estimating a nonparametric regression function, in the spirit of \cite{Korostelev1996}.

\noindent In this Note , we are interested in the problem of establishing large deviations principles of the regression operator estimate $\widehat{r}_n(\cdot)$. The results stated in the Note  deal with pointwise and uniform large deviations probabilities of  $\widehat{r}_n(\cdot)$ from $r(\cdot)$. 

\section{Results}\label{results}
Let  $F_x(h)=P[\|X_i-x\|\leq h]$, be the cumulative distribution of the real variable $W_i=\|X_i-x\|$. As in \cite{FerratyMasVieu2007}, let $\varphi(\cdot)$ be the real valued function defined by
\begin{equation}\label{varphi}
    \varphi(u)=\mathds{E}\big\{r(X)-r(x)\big|\|X-x\|=u\big\}.
\end{equation}
Before stating our results, we will consider the following conditions.
\begin{enumerate}
  \item[(C.1)]The kernel $K$ is of class $\mathcal{C}^1$ with continuous derivative on the compact support $[0,1]$.
  \item[(C.2)]The operator $r$ verifies the following Lipschitz property:
  \begin{equation}\label{LR}
   \textrm{ There exists }\beta \textrm{ such that }\forall (u,v)\in\mathcal{X}^2,\,\exists C,\,|r(u)-r(v)|\leq C \|u-v\|^{\beta}
  \end{equation}
  \item[(C.3)]There exists three functions $\ell(\cdot)$, $\phi(\cdot)$ (supposed increasing and strictly positive and tending to zero as $h$ goes to zero) and $\zeta_0(\cdot)$ such that
      \begin{enumerate}
        \item[(i)]$F_x(h)=\ell(x)\phi(h)+o(\phi(h)),$
        \item[(ii)]for all $u\in[0,1]$,  $\displaystyle{\lim_{h\rightarrow0}\frac{\phi(uh)}{\phi(h)}=:\lim_{h\rightarrow0}\zeta_h(u)=\zeta_0(u)}$.
      \end{enumerate}
  \item[C.4] $\varphi'(0)$ exists.
\end{enumerate}
\begin{remark}
There exists many examples fulfilling the decomposition mentioned in condition ({\rm{C.3}}), see for instance Proposition 1 in \cite{FerratyMasVieu2007}. The conditions stated above are classical in nonparametric estimation for functional data, we refer to \cite{FerratyMasVieu2007} and the references therein. However this condition is somewhat restrictive but is necessary to prove our results, we refer among others to \cite{FerratyLaksaciTadjVieu2010} for less restrictive conditions.
\end{remark}
\noindent Let now introduce the following functions,
\begin{equation}\label{I}
I(t)=t\lambda  \int_0^1K'(u)\exp\{ -t\lambda K(u)\}\zeta_0(u)~{\rm{d}}u;
\end{equation}
$\displaystyle{\Gamma_x^+(\lambda)=\inf_{t>0}\{\ell(x)I(t)\}}$; $\displaystyle{\Gamma_x^-(\lambda)=\inf_{t>0}\{\ell(x)I(-t)\}}$ and $\displaystyle{\Gamma_x(\lambda)=\max\{\Gamma_x^+(\lambda);\Gamma_x^-(\lambda)\}}$.

\noindent Let $x$ be a an element of the functional space $\mathcal{X}$ and $\lambda>0$. Our first theorem deals with pointwise large deviations probabilities.
\begin{theorem}\label{theorem1}
Assume that the conditions {\rm{(C.1)--(C.4)}} are satisfied. If $n\phi(h)\longrightarrow\infty$, then for any $\lambda>0$ and any $x\in\mathcal{X}$, we have
\begin{enumerate}
  \item[{\rm{(a)}}]
\begin{equation}\label{eq1Th1}
\lim_{n\rightarrow\infty}\frac{1}{n\phi(h)}\log P\big(\widehat{r}_n(x)-r(x)>\lambda\big)=\Gamma_x^+(\lambda),
\end{equation}

\item[{\rm{(b)}}]\begin{equation}\label{eq2Th1}
\lim_{n\rightarrow\infty}\frac{1}{n\phi(h)}\log P\big(\widehat{r}_n(x)-r(x)<-\lambda\big)=\Gamma_x^-(\lambda),
\end{equation}
\item[{\rm{(c)}}]\begin{equation}\label{eq3Th1}
\lim_{n\rightarrow\infty}\frac{1}{n\phi(h)}\log P\big(|\widehat{r}_n(x)-r(x)|>\lambda\big)=\Gamma_x(\lambda).
\end{equation}
\end{enumerate}
\end{theorem}
To establish uniform large deviations principles for the regression estimator we need the following additional assumptions.
Let $\mathcal{C}$ be some compact subset of $\mathcal{X}$ and $B(z_k,\xi)$ a ball centered at $z_k\in\mathcal{X}$ with radius $\xi$,  such that
for any $\xi>0$,
\begin{subequations}\label{Topo}
\begin{equation}
        \mathcal{C}\subset\bigcup_{k=1}^{\tau}B(z_k,\xi),
        \label{subeqnpart}
\end{equation}
\begin{equation}
        \exists \alpha>0,\quad  \exists C>0,\quad \tau\xi^{\alpha}=C.
        \label{geometric}
\end{equation}
\end{subequations}

\noindent Before stating the Theorem about the uniform version of our result, we introduce the following function
\begin{equation}\label{g}
g(\lambda)=\sup_{x\in\mathcal{C}}\Gamma_x(\lambda).
\end{equation}
\begin{theorem}\label{theorem2}
Assume that the conditions {\rm{(C.1)--(C.4)}} are satisfied. If $n\phi(h)\longrightarrow\infty$, then for any compact set $\mathcal{C}\subset\mathcal{X}$ satisfying conditions (\ref{Topo}) and for any $\lambda>0$,
\begin{equation}\label{eq1Th2}
\lim_{n\rightarrow\infty}\frac{1}{n\phi(h)}\log P\big(\sup_{x\in\mathcal{C}}|\widehat{r}_n(x)-r(x)|>\lambda\big)=g(\lambda)
\end{equation}
\end{theorem}
\begin{remark}
\begin{itemize}
\item[(i)] The above conditions on the covering of the compact set $\mathcal{C}$ by a finite number of balls, the geometric link between the number of balls $\tau$ and the radius $\xi$ are necessary to prove uniform convergence in the context of functional non-parametric regression and many functional non-parametric settings, see the discussion in \cite{FerratyVieu2008}.
\item[(ii)] More recently, \cite{FerratyLaksaciTadjVieu2010}  proved a very nice result about the uniform almost complete convergence of $\widehat{r}_n(\cdot)$ by use of the Kolmogorov's $\epsilon$-entropy of $\mathcal{C}$. It will be of interest to prove Theorem \ref{theorem2}, under their conditions.
\end{itemize}
\end{remark}
\section{Proofs}
\subsection{Proof of Theorem \ref{theorem1}}
We only prove the statement (\ref{eq1Th1}), (\ref{eq2Th1}) is derived in the same way.

\noindent{\rm{(a)}} Write
$$Z_{n}= \sum_{i=1}^n \big\{Y_i-r(x)-\lambda\big\}K_h(\|X_i-x\|)$$
Define $\Phi_x^{(n)}(t):=\mathds{E}\exp(tZ_{n})$ to be the moment generating function of $Z_{n}$. To prove the large deviations principles, we seek the limit of $\displaystyle{\frac{1}{n\phi(h)}\log\Phi_x^{(n)}(t)}$ as $n\longrightarrow\infty$.

\noindent Observe that
\begin{equation*}
\Phi_x^{(n)}(t)=\bigg\{1+\mathds{E}\bigg(\exp\{t[r(X_1)-r(x)-\lambda]K_h(\|X_1-x\|)\}-1\bigg)\bigg\}^n.
\end{equation*}
Using the definition of the function $\varphi$ in (\ref{varphi}), we can write
\begin{eqnarray*}
\Phi_x^{(n)}(t)&=&\bigg\{1+\mathds{E}\bigg(\exp\{t[\varphi(\|X_1-x\|)-\lambda]K_h(\|X_1-x\|)\}-1\bigg)\bigg\}^n,\\
&=&\bigg\{1+\int_0^h\bigg(\exp\{ t[\varphi(u)-\lambda]K_h(u)\}-1\bigg)~{\rm{d}}F_x(u)\bigg\}^n.\\
&=&\bigg\{1+\int_0^1\bigg(\exp\{ t[\varphi(hu)-\lambda]K(u)\}-1\bigg)~{\rm{d}}F_x(hu)\bigg\}^n.\\
\end{eqnarray*}
By {\rm{(C.4)}}, using a first order Taylor expansion of $\varphi$ about zero, we obtain
\begin{eqnarray*}
\Phi_x^{(n)}(t)&=&\bigg\{1+\int_0^1\bigg(\exp\{ t[hu\varphi'(0)-\lambda+o(1)]K(u)\}-1\bigg)~{\rm{d}}F_x(hu)\bigg\}^n,\\
\end{eqnarray*}
Integrating by parts and by {\rm{(C.1)}}, we have
\begin{eqnarray*}
\Phi_x^{(n)}(t)&=&\bigg\{1-\int_0^1t[h\varphi'(0)K(u)+K'(u)(uh\varphi'(0)-\lambda)]\\
&&\exp\{ t[hu\varphi'(0)-\lambda]K(u)\}F_x(hu)~{\rm{d}}u\bigg\}^n.\\
\end{eqnarray*}
Therefore,
\begin{eqnarray*}
\log\Phi_x^{(n)}(t)&=&n\log\bigg\{1-\int_0^1t[h\varphi'(0)K(u)+K'(u)(uh\varphi'(0)-\lambda)]\\
&&\exp\{ t[hu\varphi'(0)-\lambda]K(u)\}F_x(hu)~{\rm{d}}u\bigg\}.\\
\end{eqnarray*}
Using Taylor expansion of $\log(1+v)$ about $v=0$, we obtain
\begin{eqnarray*}
\log\Phi_x^{(n)}(t)&=&n\bigg\{-\int_0^1t[h\varphi'(0)K(u)+K'(u)(uh\varphi'(0)-\lambda)]\\
&&\exp\{ t[hu\varphi'(0)-\lambda]K(u)\}F_x(hu)~{\rm{d}}u+O(h)\bigg\}.\\
\end{eqnarray*}
Hence, from Assumption {\rm{(C.4)}} (ii) it follows that
\begin{eqnarray*}
\lim_{n\rightarrow\infty}\frac{1}{n\phi(h)}\log\Phi_x^{(n)}(t)&=&\ell(x)\bigg\{t\lambda \int_0^1K'(u)\exp\{ -t\lambda K(u)\}\zeta_0(u)~{\rm{d}}u\bigg\}\\
&=:&\ell(x)I(t).
\end{eqnarray*}
Using Theorem of \cite{PlachkySteinebach1975}, the proof of the theorem can be completed as in \cite{Louani1998}.

\noindent {\rm{(c)}} Observe that for any $x\in\mathcal{X}$,
\begin{equation*}
\max\{P(\widehat{r}_n(x)-r(x)>\lambda);\,P(\widehat{r}_n(x)-r(x)<-\lambda)\}\leq P(|\widehat{r}_n(x)-r(x)|>\lambda)
\end{equation*}
and
\begin{equation*}
P(|\widehat{r}_n(x)-r(x)|>\lambda)\leq 2\max\{P(\widehat{r}_n(x)-r(x)>\lambda);\,P(\widehat{r}_n(x)-r(x)<-\lambda)\}.
\end{equation*}
Hence,
\begin{equation*}
\Gamma_x(\lambda)\leq \lim_{n\rightarrow\infty}\frac{1}{n\phi(h)}\log P(|\widehat{r}_n(x)-r(x)|>\lambda)\leq \max\{\Gamma_x^{+}(\lambda);\,\Gamma_x^{-}(\lambda)\}=\Gamma_x(\lambda),
\end{equation*}
Thus the proof is complete.
\subsection{Proof of Theorem \ref{theorem2}}
First for any $x_0\in\mathcal{C}$, by Theorem \ref{theorem1} we have
\begin{eqnarray*}
\liminf_{n\rightarrow\infty}\frac{1}{n\phi(h)}\log P\big(\sup_{x\in\mathcal{C}}|\widehat{r}_n(x)-r(x)|>\lambda\big)&\geq & \liminf_{n\rightarrow\infty}\frac{1}{n\phi(h)}\log P\big(|\widehat{r}_n(x_0)-r(x_0)|>\lambda\big)\\
&\geq&  \Gamma_{x_0}(\lambda).
\end{eqnarray*}
Hence
\begin{equation}\label{borneinf}
\liminf_{n\rightarrow\infty}\frac{1}{n\phi(h)}\log P\big(\sup_{x\in\mathcal{C}}|\widehat{r}_n(x)-r(x)|>\lambda\big)\geq g(\lambda).
\end{equation}
To prove the reverse inequality, we note that by conditions (\ref{Topo}) it follows
\begin{equation}\label{sup1}
  \sup_{x\in\mathcal{C}}|\widehat{r}_n(x)-r(x)|\leq\max_{1\leq k \leq \tau}\sup_{x\in B(z_k,\xi)}|\widehat{r}_n(x)-r(x)|.
\end{equation}
Hence,
\begin{equation}\label{sup2}
\sup_{x\in B(z_k,\xi)}|\widehat{r}_n(x)-r(x)|\leq \sup_{x\in B(z_k,\xi)}|\widehat{r}_n(x)-\widehat{r}_n(z_k)|+ \sup_{x\in B(z_k,\xi)}|r(z_k)-r(x)|+ |\widehat{r}_n(z_k)-r(z_k)|.
\end{equation}
Using the fact that
$K$ is Lipschitz by condition {\rm{(C.1)}}, there exists $C>0$ so that
$$\sup_{x\in B(z_k,\xi)}|\widehat{r}_n(x)-\widehat{r}_n(z_k)|\leq \frac{C\xi}{n\phi(h)h}\sum_{i=1}^n|Y_i|.$$
For $n$ sufficiently large, we choose $\xi$ according to the preassigned $\epsilon>0$, so that
\begin{equation}\label{sup3}
\sup_{x\in B(z_k,\xi)}|\widehat{r}_n(x)-\widehat{r}_n(z_k)|\leq \epsilon.
\end{equation}
Moreover, $r$ is Lipschitz, hence for a suitable choice of $\xi$
\begin{equation}\label{sup4}
\sup_{x\in B(z_k,\xi)}|r(z_k)-r(x)|\leq \epsilon .
\end{equation}
Finally, (\ref{sup1})-(\ref{sup4}) yield
\begin{equation}\label{sup5}
\sup_{x\in\mathcal{C}}|\widehat{r}_n(x)-r(x)|\leq\max_{1\leq k \leq \tau}\big\{ |\widehat{r}_n(z_k)-r(z_k)|+2\epsilon\big\},
\end{equation}
which implies
$$P\big(\sup_{x\in\mathcal{C}}|\widehat{r}_n(x)-r(x)|>\lambda\big)\leq\sum_{k=1}^{ \tau}P\big(|\widehat{r}_n(z_k)-r(z_k)|>\lambda-2\epsilon\big).$$
Thus,
$$P\big(\sup_{x\in\mathcal{C}}|\widehat{r}_n(x)-r(x)|>\lambda\big)\leq\tau \max_{1\leq k\leq \tau}P\big(|\widehat{r}_n(z_k)-r(z_k)|>\lambda-2\epsilon\big).$$
It follows, by Theorem \ref{theorem1}, that
$$\limsup_{n\rightarrow\infty}\frac{1}{n\phi(h)}\log P\big(\sup_{x\in\mathcal{C}}|\widehat{r}_n(x)-r(x)|>\lambda\big)\leq \inf_{t>0} \sup_{x\in\mathcal{C}}\ell(x)I_{\epsilon}(t),$$
where
$$I_{\epsilon}(t)=t(\lambda-2\epsilon)\int_0^1K'(u)\exp\{ -t(\lambda-2\epsilon) K(u)\}\zeta_0(u)~{\rm{d}}u.$$
By continuity arguments, and the fact that
$$\inf_{t>0}\sup_{x\in\mathcal{C}}\ell(x)I(t)=\sup_{x\in\mathcal{C}}\inf_{t>0}\ell(x)I(t),$$
we obtain
\begin{equation}\label{bornesup}
\limsup_{n\rightarrow\infty}\frac{1}{n\phi(h)}\log P\big(\sup_{x\in\mathcal{C}}|\widehat{r}_n(x)-r(x)|>\lambda\big)\leq g(\lambda).
\end{equation}

\noindent Combining (\ref{borneinf}) and (\ref{bornesup}), we see that the limit exists which is $g(\lambda)$.

\end{document}